%% file: sample-inoc2024.tex
\documentclass[sigconf,edbt]{acmartinoc2024}
\def\BibTeX{{\rm B\kern-.05em{\sc i\kern-.025em b}\kern-.08em
    T\kern-.1667em\lower.7ex\hbox{E}\kern-.125emX}}

\usepackage{colortbl}
\usepackage{booktabs}
\usepackage{graphicx}
\usepackage{subcaption}
\newcolumntype{B}{>{\bfseries\columncolor{gray!30}}r}

\setcopyright{inoc2024}

\acmDOI{}

\acmISBN{}

\acmConference[INOC 2024]{11th International Network Optimization Conference (INOC)}{March 11 - 13, 2024}{Dublin, Ireland} 
\acmYear{2024}

\begin{document}
\title{Utilizing Graph Sparsification for Pre-processing in Maxcut QUBO Solver}

\author{Vorapong Suppakitpaisarn}
\affiliation{%
  \institution{The University of Tokyo}
  \city{Tokyo} 
  \country{Japan}
}
\email{vorapong@is.s.u-tokyo.ac.jp}

\author{Jin-Kao Hao}
\affiliation{%
  \institution{University of Angers}
  \city{Angers} 
  \country{France}
}
\email{jin-kao.hao@univ-angers.fr}

% The default list of authors is too long for headers}
% \renewcommand{\shortauthors}{B. Trovato et al.}
\renewcommand{\shortauthors}{}

\begin{abstract}
We suggest employing graph sparsification as a pre-processing step for maxcut programs using the QUBO solver. Quantum(-inspired) algorithms are recognized for their potential efficiency in handling quadratic unconstrained binary optimization (QUBO). Given that maxcut is an NP-hard problem and can be readily expressed using QUBO, it stands out as an exemplary case to demonstrate the effectiveness of quantum(-inspired) QUBO approaches. Here, the non-zero count in the QUBO matrix corresponds to the graph's edge count. Given that many quantum(-inspired) solvers operate through cloud services, transmitting data for dense graphs can be costly. By introducing the graph sparsification method, we aim to mitigate these communication costs. Experimental results on classical and quantum-inspired solvers indicate that this approach substantially reduces communication overheads and yields an objective value close to the optimal solution.
\end{abstract}




\maketitle

\section{Introduction}
Quantum or quantum-inspired computing is considered to have the potential to enhance the efficiency of solving various computational problems \cite{shor2002introduction,herrero2017quantum,gharibian2022dequantizing}. Specifically, these types of computers are thought to offer more effective algorithms for tackling the quadratic unconstrained binary optimization (QUBO) \cite{harwood2021formulating,oshiyama2022benchmark}. Given that several combinatorial and network optimization problems can be reformulated as QUBO, numerous researchers are actively exploring the most proficient methods for addressing these optimization issues with the aid of quantum(-inspired) QUBO solvers \cite{codognet2021constraint,codognet2022domain}.

Researchers are particularly drawn to the maximum cut problem (maxcut) \cite{mirka2023experimental,king2022improved} because it is an NP-hard problem \cite{karp2010reducibility} that can be easily expressed within the QUBO framework \cite{dunning2018works,rehfeldt2023faster}. It has been observed that preprocessing the input before feeding it to QUBO solvers can yield good solutions more efficiently than using the original data directly. As a result, various studies have introduced preprocessing strategies specifically designed for the maxcut problem to enhance the solution process \cite{ferizovic2020engineering,ferizovic2019practical,lamm2022scalable}. 

Although minimizing computation time is important for solving the maxcut problem, there is an additional challenge in addressing the problem with quantum-inspired QUBO solvers. Since quantum-inspired computers will not be commercially available for the next several decades, we are compelled to utilize these solvers through cloud services. This requires us to transmit our problems to the service providers, a step which often results in communication becoming a significant bottleneck \cite{oku2020reduce,kikuchi2023dynamical}. Therefore, our focus in this paper is on diminishing the costs associated with this communication.

\subsection{Our Contributions}

We notice that the communication cost of the maxcut problem strongly related to the number of edges in the input graph. We therefore propose to use the graph sparsification technique by the effective resistance edge sampling \cite{karger1993global,benczur1996approximating,spielman2008graph} to reduce the communication cost. 

The effective resistance technique has been demonstrated to effectively reduce the number of edges in a graph while preserving the cut size. Building on this, we propose to apply graph sparsification prior to submitting the graph to QUBO solvers to achieve results that yield cuts of large size from the QUBO solvers.

Our experiments on classical and quantum-inspired solvers validate this theoretical finding. Across all networks tested, we consistently achieved solutions where the cut size was at most 10\% smaller than the maximum cut. Remarkably, this was accomplished while reducing the number of edges – and consequently, the communication cost – by as much as 90\%. 
%In the quantum solver, our initial findings indicate that sparsification can actually give a larger cut. This improvement is attributed to the fact that sparsification leads to the creation of shallower quantum circuits, which in turn results in computations with reduced noise.

We have also noticed a decrease in computation time when using classical QUBO solvers like Gurobi \cite{gurobi2021gurobi} on max cut instances where the edges have been sparsified using our sampling technique. For instance, while Gurobi could not complete the task on the original, denser max cut problem within two hours, it was able to finish in under two seconds after we had eliminated 90\% of the edges. However, we do not consider this improvement as significant, because these solvers can still quickly find a reasonably good solution for both the original and the sparsified max cut instances. The reason Gurobi does not terminate with dense input graphs is primarily due to the extensive time required to prove the optimality of its solution.

There exists an algorithm specifically designed for addressing max cut and QUBO problems on sparse graphs within classical computing environments, referred to as McSparse, detailed in \cite{CJMM22}. We believe that our preprocessing technique could improve the computation time of the McSparse algorithm, particularly when applied to denser graph structures.  

\subsection{Related Works}

The max cut problem has garnered widespread interest among researchers, leading to the development of numerous approximation and exact algorithms. Prominent among these are the well-known SDP relaxation algorithm \cite{goemans1995improved,mahajan1999derandomizing} and algorithms for specific graph types \cite{mccormick2003easy,grotschel1984polynomial,liers2012partitioning,shih1990unifying}. In this paper, however, our focus is not on the algorithms for solving the max cut problem itself, but rather on its preprocessing. Consequently, our algorithm is designed to be compatible with all these various algorithms. 

As outlined in \cite{punnen2022quadratic}, several preprocessing techniques for QUBO solvers have been developed. Among the most significant are those based on autarkies and persistencies, which enable the determination of some binary variable values in the optimal solution \cite{nemhauser1975vertex,hammer1984roof}. Additionally, there are methods that utilize the upper bound of the relaxed program to enhance solver efficiency \cite{boros1990upper,elloumi2000decomposition}, as well as approaches centered around variable fixing \cite{billionnet2007using}. These methods have been shown to yield smaller QUBO instances that can exactly solve the original problems. In contrast, our paper introduces a preprocessing technique aimed at generating approximate QUBO instances. Importantly, our approach is designed to be compatible with these existing preprocessing methods.

The graph sparsification by edge sampling technique has been introduced to give an efficient algorithm for the maximum flow problem and the sparsest cut problem \cite{khandekar2009graph}. Also, it has been used as a pre-processing of the maximum cut problem in \cite{arora2019differentially}. The goal of using the sparsification in that paper is not to reduce the communication cost as in this paper but to increase the precision of publishing the maximum cut results under differential privacy.

%\subsection{Paper Organization}

%In the following section, we delve into the max cut problem and its formulation as a QUBO problem. This section also covers the graph sparsification technique, which is based on effective resistances. Subsequently, in Section 3, we detail our proposed method. Section 4 is dedicated to presenting our experimental results. Finally, we draw our conclusions in Section 5.

\section{Preliminaries}

\subsection{Max Cut Problem}

Consider a weighted graph \( (V, E, w) \), where \( V \) represents the set of nodes in the graph. The set of edges is denoted as \( E \subseteq \{\{u,v\} : u, v \in V, u \neq v\} \), and \( w: E \rightarrow \mathbb{R} \) is the weight function assigning a weight \( w(e) \) to each edge \( e \in E \). A cut in graph \( G \) is defined as any subset \( S \subseteq V \), with the weight of the cut \( S \) being \( w_G(S) = \sum\limits_{\{u,v\} \in E: u \in S, v \notin S} w(\{u,v\}) \). The max cut problem aims to find the cut in \( G \) that has the highest weight.

\subsection{QUBO Formulation for the Max Cut Problem}

The quadratic unconstrainted binary optimization (QUBO) is a mathematical programming in the following form 
$$\max \sum_{u} \sum_{v} Q_{u,v} x_u x_v$$
subject to $x_u \in \{0, 1\}$ for all $u$.

To express the max cut problem stated in the previous section using QUBO, we let $x_u = 1$ if $u \in S$ and $x_u = 0$ otherwise. Also, let $w'(\{u,v\}) = w(\{u,v\})$ when $\{u,v\} \in E$ and $w'(\{u,v\}) = 0$ otherwise. Since $x_u^2 = x_u$ when $x_u \in \{0,1\}$, the weight of a cut $S$ is then 
\begin{eqnarray*}
w(S) & = & \sum\limits_{\{u,v\} \in E: x_u = 1, x_v = 0} w(\{u,v\}) \\
& = & \sum\limits_{\{u,v\} \in E} w(\{u,v\}) x_u (1-x_v)\\
& = & \sum\limits_{u,v} w'(\{u,v\}) x_u (1-x_v)\\
& = & \sum_u \left[ \sum_v w'(\{u,v\}) \right] x_u - \sum_{u \neq v} w'(\{u,v\}) x_u x_v\\
& = & \sum_u \left[ \sum_v w'(\{u,v\}) \right] x_u^2 - \sum_{u \neq v} w'(\{u,v\}) x_u x_v.
\end{eqnarray*}
By defining \( Q_{u,u} = \sum\limits_v w'(\{u,v\}) \) and \( Q_{u,v} = -w'(\{u,v\}) \) for \( u \neq v \), we establish that the objective value of the QUBO corresponds to the cut size, which is also the objective value of the max cut problem. Consequently, maximizing this objective value leads to an optimal solution for the max cut problem.

In the context of solving the max cut problem with QUBO solvers available through cloud services, it becomes necessary to transmit the values of \( Q_{u,v} \) for every pair of \( u, v \). Consequently, the quantity of real numbers required to be sent is on the order of \( O(|V|^2) \). This count becomes substantially large for large graphs, turning the communication aspect into a critical bottleneck for the max cut solver.

By opting to submit only the non-zero entries of \( Q_{u,v} \), we can significantly reduce the communication cost. This means sending the QUBO problem in the format \( \left(u,v,Q_{u,v}\right) \) where \( Q_{u,v} \neq 0 \). From our definition of \( Q_{u,v} \), it is evident that for \( u \neq v \), \( Q_{u,v} \) is non-zero if and only if \( \{u,v\} \in E \). Therefore, the communication cost with this method of submission is \( O(|E|) \), which is substantially more efficient for scenarios where \( |E| \ll |V|^2 \), or in other words, when the input graph is sparse.

\subsection{Graph Sparsification by Effective Resistances \cite{spielman2008graph}}

In this section, we explore the concept of graph sparsification through Effective Resistances. Consider the input graph denoted as \( G = (V, E, w) \). Our objective is to construct a graph \( \mathsf{G} = (V, \mathsf{E}, \mathsf{w}) \) in such a way that for any cut \( S \subseteq V \), the relationship \( w_G(S) \approx w_{\mathsf{G}}(S) \) holds true. This approach aims to ensure that the weight of any given cut \( S \) in the original graph \( G \) closely approximates the weight of the same cut in the sparsified graph \( \mathsf{G} \).

Given a parameter \( q \), this method begins with an initially empty set \( \mathsf{E} \). The process involves selecting edges from the graph \( G \) a total of \( q \) times to be added to \( \mathsf{E} \). During each selection, every edge \( e \in E \) has a chance of being chosen, with this probability denoted as \( p_e \) and to be detailed in the following paragraph.

If an edge \( e \) that is not already in \( \mathsf{E} \) is selected, we assign its weight in \( \mathsf{G} \) as \( \mathsf{w}(e) = w(e)/(q \cdot p_e) \). In cases where \( e \) is already in \( \mathsf{E} \), we instead increase \( \mathsf{w}(e) \) by \( w(e)/(q \cdot p_e) \). This approach ensures that each edge's contribution to the total weight is adjusted based on its probability of selection and the number of selections, thereby maintaining a balanced representation of the graph's structure in \( \mathsf{G} \).

To establish the probability distribution \( (p_e)_{e \in E} \), we start by defining the concept of effective resistance for each edge \( e \) in \( E \), denoted as \( R_e \). We treat the graph \( G \) as if it were an electrical circuit, where each edge \( e \) is equivalent to a resistor, the resistance of which is inversely proportional to the weight of the edge, given as \( 1/w_e \). In this analogy, the effective resistance \( R_e \) of an edge \( e = \{u,v\} \) is understood as the electrical resistance experienced between nodes \( u \) and \( v \).

Subsequently, the probability \( p_e \) for each edge \( e \) is defined as the ratio of its effective resistance \( R_e \) to the sum of the effective resistances of all edges in the graph, mathematically expressed as \( p_e = R_e / \sum\limits_{e' \in E} R_{e'} \). This formulation assigns higher probabilities to edges with greater effective resistance, reflecting their relative importance in the electrical flow analogy of the graph.

The following theorem is shown in \cite{spielman2008graph}. 
\begin{theorem}
    If $q = 9 |V| \cdot \log |V| / \epsilon^2$, then, for all $S \subseteq V$,
    $$w_G(S) \leq w_{\mathsf{G}}(S) \leq (1 + \epsilon) w_G(S).$$\label{thm1}
\end{theorem}
We have from the theorem that we would obtain a sparse graph with $|\mathsf{E}| = O(|V|\log |V|)$ that preserves the cut size by the sparsification technique.

\section{Proposed Method}

\begin{figure}[h]
\centering
\includegraphics[width=0.3\textwidth]{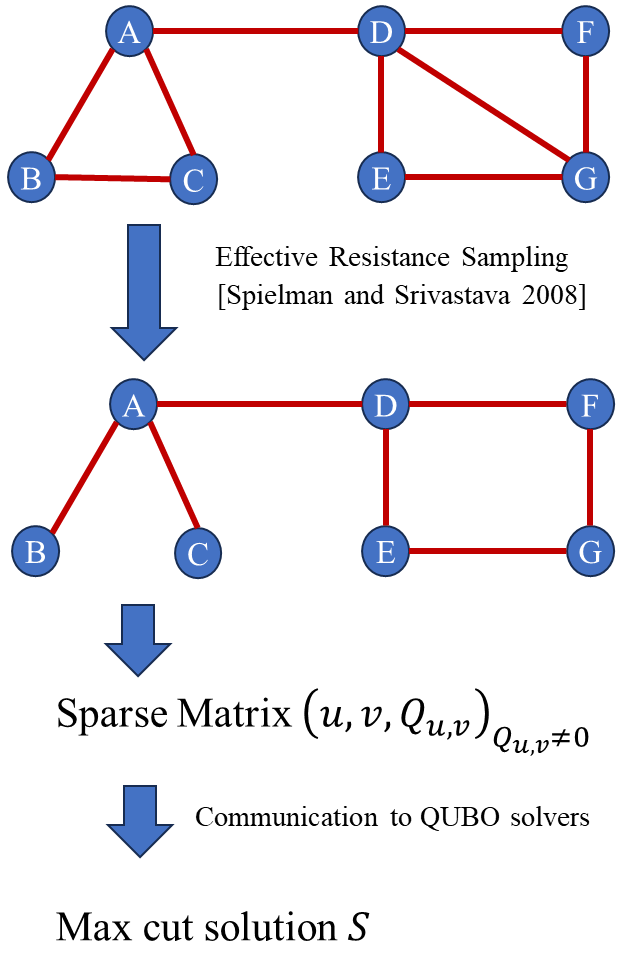}
\caption{Outline of our proposal}
\label{fig:image1}
\end{figure}

Our approach is depicted in Figure \ref{fig:image1}. Rather than directly sending the original graph \( G \) to the QUBO solver provided by cloud services, we initially apply effective resistance sampling to sparsify the graph. The resultant sparsified graph, denoted as \( \mathsf{G} \), is then submitted to the solver. 

The following theorem is directly followed from Theroem \ref{thm1}.
\begin{theorem}
    Given that \( S' \) is a cut derived from the QUBO solver using our method, and \( S^* \) represents the optimal maximum cut, we can establish that:
    $$w_G(S') \leq w_G(S^*) \leq (1 + \epsilon)w_G(S').$$
Consequently, our algorithm is a \( (1 + \epsilon) \)-approximation algorithm for the max cut problem. \label{thm2}
\end{theorem}
\begin{proof}
    Because $S'$ is the optimal maxcut solution for the graph $\mathsf{G}$, we have that $w_\mathsf{G}(S') \geq w_\mathsf{G}(S^*)$. Applying Theorem \ref{thm1}, we obtain 
    $$w_G(S') \geq \frac{1}{1 + \epsilon} w_{\mathsf{G}}(S') \geq \frac{1}{1 + \epsilon} w_{\mathsf{G}}(S^*) \geq \frac{1}{1 + \epsilon} w_{G'}(S^*).$$
    Hence, $w_G(S^*) \leq (1 + \epsilon) w_G(S')$.
\end{proof}
Theorem \ref{thm1} reveals that \( |\mathsf{E}| = O(|V|\log|V|) \), indicating that the communication cost associated with sending the sparsified graph \( \mathsf{G} \) to cloud servers is also \( O(|V|\log|V|) \). Therefore, our method can achieve an asymptotic improvement in communication costs for dense input graphs where the number of edges is on the order of \( O(|V|^2) \). However, when dealing with sparse input graphs, our approach does not yield a significant reduction in communication costs.

In Theorem \ref{thm2}, we assume that our QUBO solver is exact, meaning it always delivers the optimal solution. However, this result can be extended to scenarios where the solver is approximate. If the QUBO solver functions as an \(\alpha\)-approximation algorithm, then the outcome produced by our method can be demonstrated to be an \(\alpha(1 + \epsilon)\)-approximation.

\section{Experimental Results}

We conduct experiments on the proposed method and give the experimental results in this section. 

All experiments were carried out on a personal computer running Windows 11, equipped with an 11th Gen Intel(R) Core(TM) i7-1165G7 @2.80GHz CPU and 16 GB of RAM. The code for these experiments was written in Python. Furthermore, we utilized publicly available datasets as provided in \cite{maxcutinstances}. However, as it is assumed by the effective resistance samplings that all weights are non-negative, the weights used in our experiments are absolute values of those provided in the publicly available datasets. The values presented in this paper represent the mean of ten separate iterations. 

\subsection{Gap in Solutions Due to Graph Sparsification}

\begin{figure*}[htbp]
\centering

\begin{subfigure}[b]{0.3\textwidth} % make sure no space between 0.3 and \textwidth
  \centering
  \includegraphics[width=\textwidth]{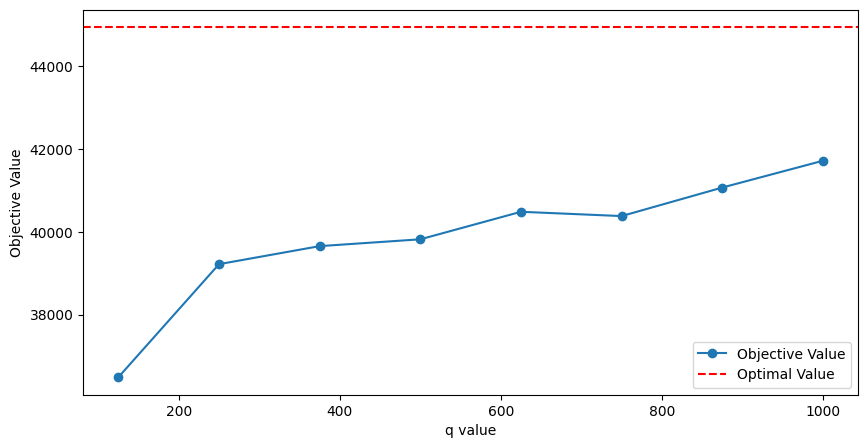}
  \caption{}
  \label{fig:2a}
\end{subfigure}
\hfill % This ensures spacing between the subfigures
\begin{subfigure}[b]{0.3\textwidth}
  \centering
  \includegraphics[width=\textwidth]{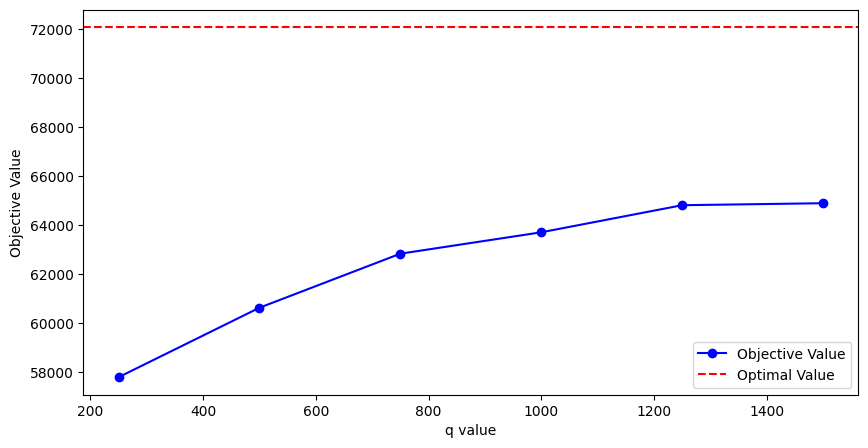}
  \caption{}
  \label{fig:2b}
\end{subfigure}
\hfill % This ensures spacing between the subfigures
\begin{subfigure}[b]{0.3\textwidth}
  \centering
  \includegraphics[width=\textwidth]{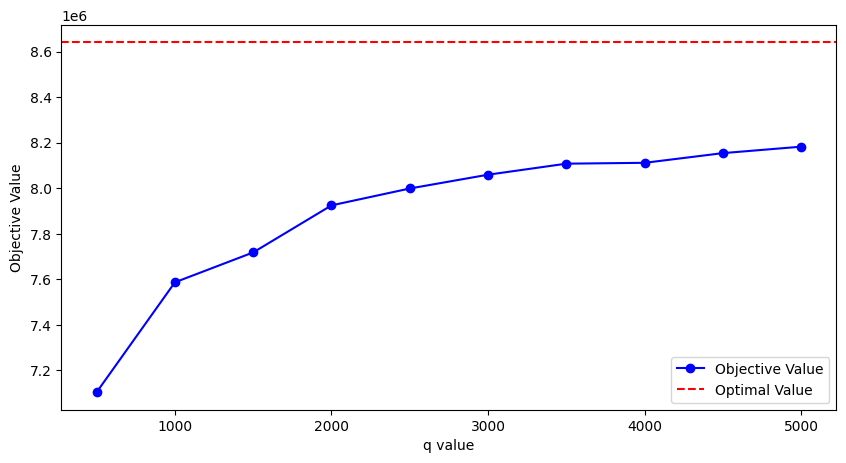}
  \caption{}
  \label{fig:2c}
\end{subfigure}

\caption{Comparisons of the optimal values derived from the original max cut instances against the objective values from the sparsified graphs for (a) be120.3.1, (b) be250.1, and (c) mannino\_k487.c}
\label{fig:2}

\end{figure*}

\begin{figure*}[htbp]
\centering

\begin{subfigure}[b]{0.3\textwidth} % make sure no space between 0.3 and \textwidth
  \centering
  \includegraphics[width=\textwidth]{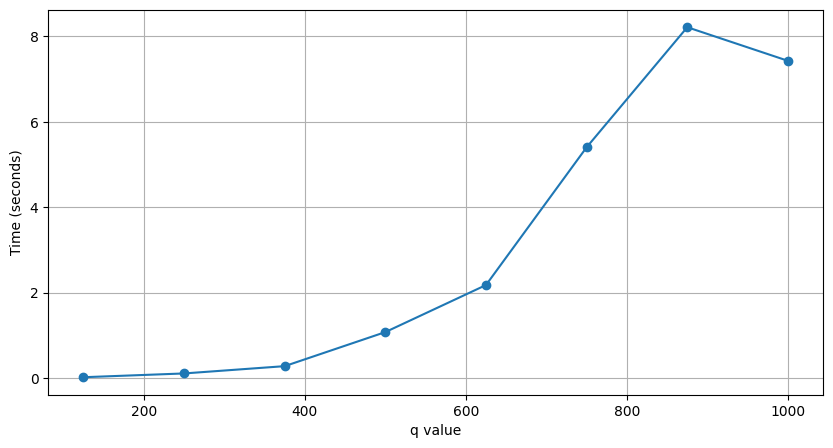}
  \caption{}
  \label{fig:3a}
\end{subfigure}
\hfill % This ensures spacing between the subfigures
\begin{subfigure}[b]{0.3\textwidth}
  \centering
  \includegraphics[width=\textwidth]{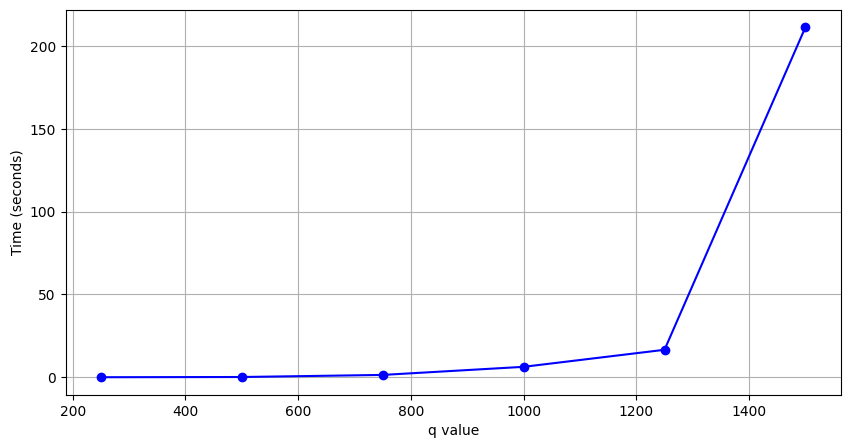}
  \caption{}
  \label{fig:3b}
\end{subfigure}
\hfill % This ensures spacing between the subfigures
\begin{subfigure}[b]{0.3\textwidth}
  \centering
  \includegraphics[width=\textwidth]{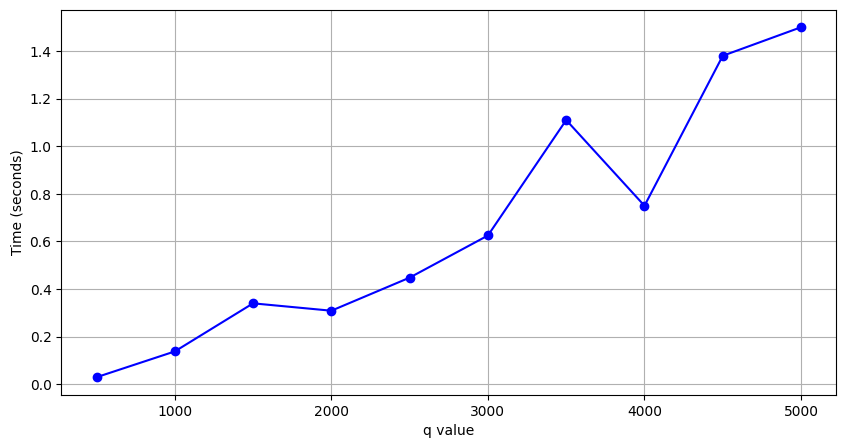}
  \caption{}
  \label{fig:3c}
\end{subfigure}

\caption{Comparisons of the computation time of Gurobi when the inputs are sparsified graphs for (a) be120.3.1, (b) be250.1, and (c) mannino\_k487.c}
\label{fig:3}

\end{figure*}

In this subsection, we examine the extent to which the optimal solutions are changed by effective resistance sampling. While our primary focus is on developing an algorithm suitable for quantum-inspired optimization, for these experiments, we have opted to use a classical solver, specifically Gurobi \cite{gurobi2021gurobi}. The rationale behind this choice is Gurobi's ability to guarantee the optimality of the solutions it generates. We employ the QUBO optimization feature available in Gurobi Optimods of Gurobi version 10.0.3.

We initially evaluated our proposed methods using three distinct instances, each varying in node count and type. These instances are "be120.3.1", "be250.1", and "mannino\_k487c". The datasets "be120.3.1" and "be250.1" are synthetic and were utilized in \cite{billionnet2007using}. They were created using generators described in \cite{pardalos1990computational}. Specifically, "be120.3.1" comprises 121 nodes and 2242 edges, whereas "be250.1" contains 251 nodes and 3269 edges. The "mannino\_k487c" dataset, on the other hand, is rooted in real-world data concerning radio frequency interferences among major Italian cities, as detailed in \cite{bonato2014lifting}. This dataset features 487 nodes and 8511 edges.

In these experiments, we focus on the variable \( q \), which represents the number of times edges are sampled from the input graph. We experiment with varying the value of \( q \). It is crucial to understand that \( q \) does not directly correspond to the number of edges in the sparsified graph, denoted as \( |\mathsf{E}| \), since an edge can be selected multiple times during sampling. However, it is evident that \( |\mathsf{E}| \leq q \), and generally, a larger \( q \) tends to result in a higher number of edges in the sparsified graph.

Theorem \ref{thm1} suggests setting \( q \) to \( \frac{9|V| \cdot \log |V|}{\epsilon^2} \). While this value is theoretically smaller than the edge count for dense graphs in an asymptotic sense, the sizeable constant factor \( \frac{9}{\epsilon^2} \) can lead to a \( q \) that exceeds the actual number of edges, especially when the input graph \( G \) is relatively small. Take, for instance, when \( \epsilon \) is 0.1, this results in approximately 522261 for "be120.3.1", 1248200 for "be250.1", and 2712316 for "mannino\_k487c". Because of this, we have opted to use a reduced \( q \) value for our experimental evaluations. 

Figure \ref{fig:2} presents the outcomes for the three specified instances. Examination of the figure reveals that the objective function improves as more edges are sampled and the value of \( q \) increases. Our method achieves a cut size that exceeds 90\% of the optimal cut size for \( q \geq 500 \) in "be120.3.1", for \( q \geq 1500 \) in "be250.1", and for \( q \geq 2000 \) in "mannino\_k487c". Correspondingly, these thresholds yield edge counts of 424 for "be120.3.1", 1092 for "be250.1", and 1231 for "mannino\_k487c". These results indicate that our approach not only secures a 0.9-approximation to the solution but also facilitates a substantial reduction in communication costs—81\% for "be120.3.1", 67\% for "be250.1", and 86\% for "mannino\_k487c". 

\begin{table*}[ht]
\centering
\caption{The reduction in communication costs and the approximation ratios achieved by our algorithms across various graph types are detailed in \cite{maxcutinstances}}
\begin{tabular}{@{}lrrrBrrB@{}}
\toprule
Dataset Name & $|V|$ & $|E|$ & $|\mathsf{E}|$ & \textbf{Reduction} & Optimal Value & Our Objective Value & \textbf{Approx. Ratio} \\ 
& & & & \textbf{in Comm. Cost} & & & \\ 
\midrule
bqp250-1 & 251 & 3339 & 1163.6 & 0.65151 & 143669 & 129863 & 0.90390 \\
gka1e & 201 & 2124 & 810.7 & 0.61831 & 48263 & 42829 & 0.88741 \\
ising2.5-150\_6666 & 150 & 10722 & 387.3 & 0.96388 & 9067341 & 8502808 & 0.93774 \\
g05\_100.0 & 100 & 2475 & 452.3 & 0.81725 & 1430 & 1309 & 0.91538 \\
w05\_100.0 & 100 & 2343 & 432.3 & 0.81549 & 7737 & 7033.9 & 0.90912 \\
G\_1 & 800 & 19176 & 3598.3 & 0.8124
 & 11598 & 10412.3 & 0.8978 \\
\bottomrule
\end{tabular}
\label{tab1}
\end{table*}

Our experiment results with these three datasets yield a 0.9-approximation solution when setting \( q \) to roughly \( 5|V| \). Consequently, we extrapolate this finding to additional instance types in~\cite{maxcutinstances}. As demonstrated in Table \ref{tab1}, a similar approximation ratio is achieved for all tested instance types with \( q \) set at \( 5|V| \). Notably, there is a substantial decrease in communication cost particularly when the original graph \( G \) is dense.

\subsection{Computation Time in Classical Solver}

Because graph sparsification techniques are often employed to reduce computation time, it is worth investigating whether our sparsification method also reduces the computational times for classical solvers.

Figure \ref{fig:3} demonstrates that sparsification does indeed have a significant effect on reducing the computational time for Gurobi. For the original "be120.3.1" and "be250.1" inputs, the solver requires more than three hours to find a solution, whereas with the sparsified graphs at \( q \approx 5|V| \), the computation times drop dramatically to 2.18 seconds and 16.6 seconds, respectively. There is also a clear pattern where larger values of \( q \) and increased edge counts correlate with longer computation times.

Despite this, it is noteworthy that Gurobi is able to quickly find reasonably good solutions for denser graphs. In every test conducted, solutions surpassing those of the sparsified graphs were obtained in under five seconds using the original graphs. Most of the time which the solver used on the dense graphs is for proving that their solutions are optimal. This leads us to conclude that reduced computation time may not be a decisive advantage of sparsification techniques.

\subsection{Experiments on Quantum-Inspired Solvers}

We employed the Fixstars Amplify Annealing Engine \cite{fixstars} to corroborate our findings with QUBO quantum-inspired solvers. For all instances, whether original or sparsified, we imposed a solver time constraint of 10 seconds. Consistent with the methodology outlined in the preceding section, we set the value of \( q \) to be \( 5|V| \) in this experiment.

\begin{table}[ht]
\centering
\caption{Reduction in communication Cost and changes in objective value by the graph sparsification technique when conducting experiments on QUBO quantum-inspired solver}
\begin{tabular}{@{}lcc@{}}
\toprule
Dataset Name & Reduction in & Changes in \\
& Comm. Cost & Objective Value \\

\midrule
be120.3.1 & 0.779 & 0.911 \\
be250.1 & 0.704 & 0.891 \\
mannino\_k487.c & 0.835 & 0.92 \\
bqp250-1 & 0.71 & 0.895 \\
gka1e & 0.65 & 0.909 \\
ising2.5-150\_6666 & 0.964 & 0.947 \\
g05\_100.0 & 0.816 & 0.908 \\
w05\_100.0 & 0.812 & 0.902 \\
G\_1 & 0.812 & 0.897 \\
\bottomrule
\end{tabular}
\label{tab2}
\end{table}

Table \ref{tab2} shows that the outcomes obtained using quantum-inspired solvers align closely with those presented in Table \ref{tab1}, confirming consistency across all datasets tested with the classical solver.

\subsection{Discussions on Results on Classical and Quantum-Inspired Solvers}

The argument could be made that similar or even superior approximation ratios to those achieved in our research might be attainable using an approximation algorithm based on semi-definite programming, as demonstrated in previous studies \cite{goemans1995improved,mahajan1999derandomizing}. This algorithm is indeed capable of providing polynomial-time approximation solutions for maxcut problems. However, a notable limitation of semi-definite programming is its computational intensity, particularly for problems involving over 100,000 nodes \cite{kim2003exact}, where local execution becomes impractical. In such scenarios, our method proves advantageous, offering a viable solution by enabling the processing of these large instances through cloud services.

In these experiments, our primary objective is to demonstrate that edge sampling can yield reasonable approximation ratios. Therefore, we confined our experimentation to smaller instances (with \( |V| \leq 800 \)) where obtaining optimal solutions is feasible. Nonetheless, given the consistent results across all tested instance sizes, we are confident that similar outcomes would be achievable with larger graphs.

\section{Conclusion and Future Works}

In our paper, we introduce the application of graph sparsification as a preprocessing step for solving the maximum cut problem in cloud-based environments. Our experimental results demonstrate that this approach, when applied to classical and quantum-inspired solvers, consistently yields solutions with an approximation ratio of about 0.9, while simultaneously achieving a significant reduction in communication costs to cloud servers, ranging between 60\% and 95\%. %Additionally, for quantum solvers, the sparsification method may lead to enhanced outcomes. This improvement is attributed to the creation of shallower quantum circuits and the decrease in computational noise, a key advantage in quantum computing scenarios.

In our future research, we plan to expand our experiments to include quantum solvers. It is understood that a sparser graph results in shallower quantum circuits, thereby reducing the noise in quantum computations. The graph sparsification technique has already been used for solving max cut for the noisy data published under differential privacy \cite{arora2019differentially}. Based on this understanding, our hypothesis is that graph sparsification not only diminishes the communication cost but also enhances the quality of solutions derived from quantum solvers.

\begin{acks}
This work is supported by JST SICORP Grant Number JPMJSC2208, Japan. The authors would like to thank
Prof. Philippe Codognet and Prof. Hiroshi Imai for several useful comments and ideas which significantly improved this paper. 
%Part of the results in this paper were obtained using an IBM Quantum computing system as part of the IBM Quantum Hub at The University of Tokyo.
\end{acks}

\bibliographystyle{ACM-Reference-Format}
\input{sample-inoc2024.bbl}

\end{document}

%% file: sample-inoc2024.bbl
%%% -*-BibTeX-*-
%%% Do NOT edit. File created by BibTeX with style
%%% ACM-Reference-Format-Journals [18-Jan-2012].

%% file: sample-inoc2024.bbl
\begin{thebibliography}{46}

%%% ====================================================================
%%% NOTE TO THE USER: you can override these defaults by providing
%%% customized versions of any of these macros before the \bibliography
%%% command.  Each of them MUST provide its own final punctuation,
%%% except for \shownote{}, \showDOI{}, and \showURL{}.  The latter two
%%% do not use final punctuation, in order to avoid confusing it with
%%% the Web address.
%%%
%%% To suppress output of a particular field, define its macro to expand
%%% to an empty string, or better, \unskip, like this:
%%%
%%% \newcommand{\showDOI}[1]{\unskip}   % LaTeX syntax
%%%
%%% \def \showDOI #1{\unskip}           % plain TeX syntax
%%%
%%% ====================================================================

\ifx \showCODEN    \undefined \def \showCODEN     #1{\unskip}     \fi
\ifx \showDOI      \undefined \def \showDOI       #1{#1}\fi
\ifx \showISBNx    \undefined \def \showISBNx     #1{\unskip}     \fi
\ifx \showISBNxiii \undefined \def \showISBNxiii  #1{\unskip}     \fi
\ifx \showISSN     \undefined \def \showISSN      #1{\unskip}     \fi
\ifx \showLCCN     \undefined \def \showLCCN      #1{\unskip}     \fi
\ifx \shownote     \undefined \def \shownote      #1{#1}          \fi
\ifx \showarticletitle \undefined \def \showarticletitle #1{#1}   \fi
\ifx \showURL      \undefined \def \showURL       {\relax}        \fi
% The following commands are used for tagged output and should be
% invisible to TeX
\providecommand\bibfield[2]{#2}
\providecommand\bibinfo[2]{#2}
\providecommand\natexlab[1]{#1}
\providecommand\showeprint[2][]{arXiv:#2}

\bibitem[\protect\citeauthoryear{Arora and Upadhyay}{Arora and Upadhyay}{2019}]%
        {arora2019differentially}
\bibfield{author}{\bibinfo{person}{Raman Arora} {and} \bibinfo{person}{Jalaj Upadhyay}.} \bibinfo{year}{2019}\natexlab{}.
\newblock \showarticletitle{On differentially private graph sparsification and applications}.
\newblock \bibinfo{journal}{\emph{Advances in neural information processing systems}}  \bibinfo{volume}{32} (\bibinfo{year}{2019}).
\newblock


\bibitem[\protect\citeauthoryear{Bencz{\'u}r and Karger}{Bencz{\'u}r and Karger}{1996}]%
        {benczur1996approximating}
\bibfield{author}{\bibinfo{person}{Andr{\'a}s~A Bencz{\'u}r} {and} \bibinfo{person}{David~R Karger}.} \bibinfo{year}{1996}\natexlab{}.
\newblock \showarticletitle{Approximating {$s$-$t$} minimum cuts in {$O(n^2)$} time}. In \bibinfo{booktitle}{\emph{Proceedings of the twenty-eighth annual ACM symposium on Theory of computing}}. \bibinfo{pages}{47--55}.
\newblock


\bibitem[\protect\citeauthoryear{Billionnet and Elloumi}{Billionnet and Elloumi}{2007}]%
        {billionnet2007using}
\bibfield{author}{\bibinfo{person}{Alain Billionnet} {and} \bibinfo{person}{Sourour Elloumi}.} \bibinfo{year}{2007}\natexlab{}.
\newblock \showarticletitle{Using a mixed integer quadratic programming solver for the unconstrained quadratic 0-1 problem}.
\newblock \bibinfo{journal}{\emph{Mathematical programming}}  \bibinfo{volume}{109} (\bibinfo{year}{2007}), \bibinfo{pages}{55--68}.
\newblock


\bibitem[\protect\citeauthoryear{Bonato, J{\"u}nger, Reinelt, and Rinaldi}{Bonato et~al\mbox{.}}{2014}]%
        {bonato2014lifting}
\bibfield{author}{\bibinfo{person}{Thorsten Bonato}, \bibinfo{person}{Michael J{\"u}nger}, \bibinfo{person}{Gerhard Reinelt}, {and} \bibinfo{person}{Giovanni Rinaldi}.} \bibinfo{year}{2014}\natexlab{}.
\newblock \showarticletitle{Lifting and separation procedures for the cut polytope}.
\newblock \bibinfo{journal}{\emph{Mathematical Programming}}  \bibinfo{volume}{146} (\bibinfo{year}{2014}), \bibinfo{pages}{351--378}.
\newblock


\bibitem[\protect\citeauthoryear{Boros, Crama, and Hammer}{Boros et~al\mbox{.}}{1990}]%
        {boros1990upper}
\bibfield{author}{\bibinfo{person}{Emil Boros}, \bibinfo{person}{Yves Crama}, {and} \bibinfo{person}{Peter~L Hammer}.} \bibinfo{year}{1990}\natexlab{}.
\newblock \showarticletitle{Upper-bounds for quadratic 0--1 maximization}.
\newblock \bibinfo{journal}{\emph{Operations Research Letters}} \bibinfo{volume}{9}, \bibinfo{number}{2} (\bibinfo{year}{1990}), \bibinfo{pages}{73--79}.
\newblock


\bibitem[\protect\citeauthoryear{Charfreitag, J{\"u}nger, Mallach, and Mutzel}{Charfreitag et~al\mbox{.}}{2022}]%
        {CJMM22}
\bibfield{author}{\bibinfo{person}{Jonas Charfreitag}, \bibinfo{person}{Michael J{\"u}nger}, \bibinfo{person}{Sven Mallach}, {and} \bibinfo{person}{Petra Mutzel}.} \bibinfo{year}{2022}\natexlab{}.
\newblock \showarticletitle{{M}c{S}parse: {E}xact Solutions of Sparse Maximum Cut and Sparse Unconstrained Binary Quadratic Optimization Problems}. In \bibinfo{booktitle}{\emph{ALENEX 2022}}. \bibinfo{pages}{54--66}.
\newblock


\bibitem[\protect\citeauthoryear{Codognet}{Codognet}{2021}]%
        {codognet2021constraint}
\bibfield{author}{\bibinfo{person}{Philippe Codognet}.} \bibinfo{year}{2021}\natexlab{}.
\newblock \showarticletitle{Constraint solving by quantum annealing}. In \bibinfo{booktitle}{\emph{50th International Conference on Parallel Processing Workshop}}. \bibinfo{pages}{1--10}.
\newblock


\bibitem[\protect\citeauthoryear{Codognet}{Codognet}{2022}]%
        {codognet2022domain}
\bibfield{author}{\bibinfo{person}{Philippe Codognet}.} \bibinfo{year}{2022}\natexlab{}.
\newblock \showarticletitle{Domain-wall/unary encoding in QUBO for permutation problems}. In \bibinfo{booktitle}{\emph{2022 IEEE International Conference on Quantum Computing and Engineering (QCE)}}. IEEE, \bibinfo{pages}{167--173}.
\newblock


\bibitem[\protect\citeauthoryear{Dial}{Dial}{2022}]%
        {Eagle}
\bibfield{author}{\bibinfo{person}{Oliver Dial}.} \bibinfo{year}{2022}\natexlab{}.
\newblock \bibinfo{title}{Eagle’s quantum performance progress}.
\newblock
\newblock
\urldef\tempurl%
\url{https://research.ibm.com/blog/eagle-quantum-processor-performance}
\showURL{%
\tempurl}


\bibitem[\protect\citeauthoryear{Dunning, Gupta, and Silberholz}{Dunning et~al\mbox{.}}{2018}]%
        {dunning2018works}
\bibfield{author}{\bibinfo{person}{Iain Dunning}, \bibinfo{person}{Swati Gupta}, {and} \bibinfo{person}{John Silberholz}.} \bibinfo{year}{2018}\natexlab{}.
\newblock \showarticletitle{What works best when? A systematic evaluation of heuristics for Max-Cut and QUBO}.
\newblock \bibinfo{journal}{\emph{INFORMS Journal on Computing}} \bibinfo{volume}{30}, \bibinfo{number}{3} (\bibinfo{year}{2018}), \bibinfo{pages}{608--624}.
\newblock


\bibitem[\protect\citeauthoryear{Elloumi, Faye, and Soutif}{Elloumi et~al\mbox{.}}{2000}]%
        {elloumi2000decomposition}
\bibfield{author}{\bibinfo{person}{Sourour Elloumi}, \bibinfo{person}{Alain Faye}, {and} \bibinfo{person}{Eric Soutif}.} \bibinfo{year}{2000}\natexlab{}.
\newblock \showarticletitle{Decomposition and linearization for 0-1 quadratic programming}.
\newblock \bibinfo{journal}{\emph{Annals of Operations Research}} \bibinfo{volume}{99}, \bibinfo{number}{1-4} (\bibinfo{year}{2000}), \bibinfo{pages}{79--93}.
\newblock


\bibitem[\protect\citeauthoryear{Ferizovic}{Ferizovic}{2019}]%
        {ferizovic2019practical}
\bibfield{author}{\bibinfo{person}{Damir Ferizovic}.} \bibinfo{year}{2019}\natexlab{}.
\newblock \emph{\bibinfo{title}{A Practical Analysis of Kernelization Techniques for the Maximum Cut Problem}}.
\newblock \bibinfo{thesistype}{Ph.D. Dissertation}. \bibinfo{school}{Karlsruher Institut f{\"u}r Technologie (KIT)}.
\newblock


\bibitem[\protect\citeauthoryear{Ferizovic, Hespe, Lamm, Mnich, Schulz, and Strash}{Ferizovic et~al\mbox{.}}{2020}]%
        {ferizovic2020engineering}
\bibfield{author}{\bibinfo{person}{Damir Ferizovic}, \bibinfo{person}{Demian Hespe}, \bibinfo{person}{Sebastian Lamm}, \bibinfo{person}{Matthias Mnich}, \bibinfo{person}{Christian Schulz}, {and} \bibinfo{person}{Darren Strash}.} \bibinfo{year}{2020}\natexlab{}.
\newblock \showarticletitle{Engineering kernelization for maximum cut}. In \bibinfo{booktitle}{\emph{ALENEX 2020}}. SIAM, \bibinfo{pages}{27--41}.
\newblock


\bibitem[\protect\citeauthoryear{Fixstars}{Fixstars}{2023}]%
        {fixstars}
\bibfield{author}{\bibinfo{person}{Fixstars}.} \bibinfo{year}{2023}\natexlab{}.
\newblock \bibinfo{title}{About Amplify AE}.
\newblock
\newblock
\urldef\tempurl%
\url{https://amplify.fixstars.com/ja/docs/amplify-ae/about.html}
\showURL{%
\tempurl}


\bibitem[\protect\citeauthoryear{Gharibian and Le~Gall}{Gharibian and Le~Gall}{2022}]%
        {gharibian2022dequantizing}
\bibfield{author}{\bibinfo{person}{Sevag Gharibian} {and} \bibinfo{person}{Fran{\c{c}}ois Le~Gall}.} \bibinfo{year}{2022}\natexlab{}.
\newblock \showarticletitle{Dequantizing the quantum singular value transformation: hardness and applications to quantum chemistry and the quantum PCP conjecture}. In \bibinfo{booktitle}{\emph{Proceedings of the 54th Annual ACM SIGACT Symposium on Theory of Computing}}. \bibinfo{pages}{19--32}.
\newblock


\bibitem[\protect\citeauthoryear{Goemans and Williamson}{Goemans and Williamson}{1995}]%
        {goemans1995improved}
\bibfield{author}{\bibinfo{person}{Michel~X Goemans} {and} \bibinfo{person}{David~P Williamson}.} \bibinfo{year}{1995}\natexlab{}.
\newblock \showarticletitle{Improved approximation algorithms for maximum cut and satisfiability problems using semidefinite programming}.
\newblock \bibinfo{journal}{\emph{Journal of the ACM (JACM)}} \bibinfo{volume}{42}, \bibinfo{number}{6} (\bibinfo{year}{1995}), \bibinfo{pages}{1115--1145}.
\newblock


\bibitem[\protect\citeauthoryear{Gr{\"o}tschel and Nemhauser}{Gr{\"o}tschel and Nemhauser}{1984}]%
        {grotschel1984polynomial}
\bibfield{author}{\bibinfo{person}{Martin Gr{\"o}tschel} {and} \bibinfo{person}{George~L Nemhauser}.} \bibinfo{year}{1984}\natexlab{}.
\newblock \showarticletitle{A polynomial algorithm for the max-cut problem on graphs without long odd cycles}.
\newblock \bibinfo{journal}{\emph{Mathematical Programming}} \bibinfo{volume}{29}, \bibinfo{number}{1} (\bibinfo{year}{1984}), \bibinfo{pages}{28--40}.
\newblock


\bibitem[\protect\citeauthoryear{Guerreschi and Matsuura}{Guerreschi and Matsuura}{2019}]%
        {guerreschi2019qaoa}
\bibfield{author}{\bibinfo{person}{Gian~Giacomo Guerreschi} {and} \bibinfo{person}{Anne~Y Matsuura}.} \bibinfo{year}{2019}\natexlab{}.
\newblock \showarticletitle{QAOA for Max-Cut requires hundreds of qubits for quantum speed-up}.
\newblock \bibinfo{journal}{\emph{Scientific reports}} \bibinfo{volume}{9}, \bibinfo{number}{1} (\bibinfo{year}{2019}), \bibinfo{pages}{6903}.
\newblock


\bibitem[\protect\citeauthoryear{Gurobi~Optimization}{Gurobi~Optimization}{2021}]%
        {gurobi2021gurobi}
\bibfield{author}{\bibinfo{person}{LLC Gurobi~Optimization}.} \bibinfo{year}{2021}\natexlab{}.
\newblock \bibinfo{title}{Gurobi optimizer reference manual}.
\newblock
\newblock


\bibitem[\protect\citeauthoryear{Hammer, Hansen, and Simeone}{Hammer et~al\mbox{.}}{1984}]%
        {hammer1984roof}
\bibfield{author}{\bibinfo{person}{Peter~L Hammer}, \bibinfo{person}{Pierre Hansen}, {and} \bibinfo{person}{Bruno Simeone}.} \bibinfo{year}{1984}\natexlab{}.
\newblock \showarticletitle{Roof duality, complementation and persistency in quadratic 0--1 optimization}.
\newblock \bibinfo{journal}{\emph{Mathematical programming}}  \bibinfo{volume}{28} (\bibinfo{year}{1984}), \bibinfo{pages}{121--155}.
\newblock


\bibitem[\protect\citeauthoryear{Harwood, Gambella, Trenev, Simonetto, Bernal, and Greenberg}{Harwood et~al\mbox{.}}{2021}]%
        {harwood2021formulating}
\bibfield{author}{\bibinfo{person}{Stuart Harwood}, \bibinfo{person}{Claudio Gambella}, \bibinfo{person}{Dimitar Trenev}, \bibinfo{person}{Andrea Simonetto}, \bibinfo{person}{David Bernal}, {and} \bibinfo{person}{Donny Greenberg}.} \bibinfo{year}{2021}\natexlab{}.
\newblock \showarticletitle{Formulating and solving routing problems on quantum computers}.
\newblock \bibinfo{journal}{\emph{IEEE Transactions on Quantum Engineering}}  \bibinfo{volume}{2} (\bibinfo{year}{2021}), \bibinfo{pages}{1--17}.
\newblock


\bibitem[\protect\citeauthoryear{Herrero-Collantes and Garcia-Escartin}{Herrero-Collantes and Garcia-Escartin}{2017}]%
        {herrero2017quantum}
\bibfield{author}{\bibinfo{person}{Miguel Herrero-Collantes} {and} \bibinfo{person}{Juan~Carlos Garcia-Escartin}.} \bibinfo{year}{2017}\natexlab{}.
\newblock \showarticletitle{Quantum random number generators}.
\newblock \bibinfo{journal}{\emph{Reviews of Modern Physics}} \bibinfo{volume}{89}, \bibinfo{number}{1} (\bibinfo{year}{2017}), \bibinfo{pages}{015004}.
\newblock


\bibitem[\protect\citeauthoryear{Japan}{Japan}{2023}]%
        {kawasaki}
\bibfield{author}{\bibinfo{person}{IBM Japan}.} \bibinfo{year}{2023}\natexlab{}.
\newblock \bibinfo{title}{IBM Quantum System One {``ibm\_kawasaki''}}.
\newblock
\newblock
\urldef\tempurl%
\url{https://www.ibm.com/blogs/think/jp-ja/tag/ibm_kawasaki/}
\showURL{%
\tempurl}


\bibitem[\protect\citeauthoryear{Karger}{Karger}{1993}]%
        {karger1993global}
\bibfield{author}{\bibinfo{person}{David~R Karger}.} \bibinfo{year}{1993}\natexlab{}.
\newblock \showarticletitle{Global min-cuts in RNC, and other ramifications of a simple min-cut algorithm}. In \bibinfo{booktitle}{\emph{Proceedings of the fourth annual ACM-SIAM symposium on Discrete algorithms}}. \bibinfo{pages}{21--30}.
\newblock


\bibitem[\protect\citeauthoryear{Karp}{Karp}{2010}]%
        {karp2010reducibility}
\bibfield{author}{\bibinfo{person}{Richard~M Karp}.} \bibinfo{year}{2010}\natexlab{}.
\newblock \bibinfo{booktitle}{\emph{Reducibility among combinatorial problems}}.
\newblock \bibinfo{publisher}{Springer}.
\newblock


\bibitem[\protect\citeauthoryear{Khandekar, Rao, and Vazirani}{Khandekar et~al\mbox{.}}{2009}]%
        {khandekar2009graph}
\bibfield{author}{\bibinfo{person}{Rohit Khandekar}, \bibinfo{person}{Satish Rao}, {and} \bibinfo{person}{Umesh Vazirani}.} \bibinfo{year}{2009}\natexlab{}.
\newblock \showarticletitle{Graph partitioning using single commodity flows}.
\newblock \bibinfo{journal}{\emph{Journal of the ACM (JACM)}} \bibinfo{volume}{56}, \bibinfo{number}{4} (\bibinfo{year}{2009}), \bibinfo{pages}{1--15}.
\newblock


\bibitem[\protect\citeauthoryear{Kikuchi, Togawa, and Tanaka}{Kikuchi et~al\mbox{.}}{2023}]%
        {kikuchi2023dynamical}
\bibfield{author}{\bibinfo{person}{Shuta Kikuchi}, \bibinfo{person}{Nozomu Togawa}, {and} \bibinfo{person}{Shu Tanaka}.} \bibinfo{year}{2023}\natexlab{}.
\newblock \showarticletitle{Dynamical process of a bit-width reduced Ising model with simulated annealing}.
\newblock \bibinfo{journal}{\emph{arXiv preprint arXiv:2304.12796}} (\bibinfo{year}{2023}).
\newblock


\bibitem[\protect\citeauthoryear{Kim and Kojima}{Kim and Kojima}{2003}]%
        {kim2003exact}
\bibfield{author}{\bibinfo{person}{Sunyoung Kim} {and} \bibinfo{person}{Masakazu Kojima}.} \bibinfo{year}{2003}\natexlab{}.
\newblock \showarticletitle{Exact solutions of some nonconvex quadratic optimization problems via SDP and SOCP relaxations}.
\newblock \bibinfo{journal}{\emph{Computational optimization and applications}}  \bibinfo{volume}{26} (\bibinfo{year}{2003}), \bibinfo{pages}{143--154}.
\newblock


\bibitem[\protect\citeauthoryear{King}{King}{2022}]%
        {king2022improved}
\bibfield{author}{\bibinfo{person}{Robbie King}.} \bibinfo{year}{2022}\natexlab{}.
\newblock \showarticletitle{An improved approximation algorithm for quantum max-cut}.
\newblock \bibinfo{journal}{\emph{arXiv preprint arXiv:2209.02589}} (\bibinfo{year}{2022}).
\newblock


\bibitem[\protect\citeauthoryear{Kullback and Leibler}{Kullback and Leibler}{1951}]%
        {kullback1951information}
\bibfield{author}{\bibinfo{person}{Solomon Kullback} {and} \bibinfo{person}{Richard~A Leibler}.} \bibinfo{year}{1951}\natexlab{}.
\newblock \showarticletitle{On information and sufficiency}.
\newblock \bibinfo{journal}{\emph{The annals of mathematical statistics}} \bibinfo{volume}{22}, \bibinfo{number}{1} (\bibinfo{year}{1951}), \bibinfo{pages}{79--86}.
\newblock


\bibitem[\protect\citeauthoryear{Lamm}{Lamm}{2022}]%
        {lamm2022scalable}
\bibfield{author}{\bibinfo{person}{Sebastian Lamm}.} \bibinfo{year}{2022}\natexlab{}.
\newblock \emph{\bibinfo{title}{Scalable Graph Algorithms using Practically Efficient Data Reductions}}.
\newblock \bibinfo{thesistype}{Ph.D. Dissertation}. \bibinfo{school}{Karlsruher Institut f{\"u}r Technologie (KIT)}.
\newblock


\bibitem[\protect\citeauthoryear{Liers and Pardella}{Liers and Pardella}{2012}]%
        {liers2012partitioning}
\bibfield{author}{\bibinfo{person}{Frauke Liers} {and} \bibinfo{person}{Gregor Pardella}.} \bibinfo{year}{2012}\natexlab{}.
\newblock \showarticletitle{Partitioning planar graphs: a fast combinatorial approach for max-cut}.
\newblock \bibinfo{journal}{\emph{Computational Optimization and Applications}} \bibinfo{volume}{51}, \bibinfo{number}{1} (\bibinfo{year}{2012}), \bibinfo{pages}{323--344}.
\newblock


\bibitem[\protect\citeauthoryear{Mahajan and Ramesh}{Mahajan and Ramesh}{1999}]%
        {mahajan1999derandomizing}
\bibfield{author}{\bibinfo{person}{Sanjeev Mahajan} {and} \bibinfo{person}{Hariharan Ramesh}.} \bibinfo{year}{1999}\natexlab{}.
\newblock \showarticletitle{Derandomizing approximation algorithms based on semidefinite programming}.
\newblock \bibinfo{journal}{\emph{SIAM J. Comput.}} \bibinfo{volume}{28}, \bibinfo{number}{5} (\bibinfo{year}{1999}), \bibinfo{pages}{1641--1663}.
\newblock


\bibitem[\protect\citeauthoryear{Mallach, Junger, Charfreitag, and Jordan}{Mallach et~al\mbox{.}}{2021}]%
        {maxcutinstances}
\bibfield{author}{\bibinfo{person}{Sven Mallach}, \bibinfo{person}{Michael Junger}, \bibinfo{person}{Jonas Charfreitag}, {and} \bibinfo{person}{Claude Jordan}.} \bibinfo{year}{2021}\natexlab{}.
\newblock \bibinfo{title}{(Prototype of a) MaxCut and BQP Instance Library}.
\newblock
\newblock


\bibitem[\protect\citeauthoryear{McCormick, Rao, and Rinaldi}{McCormick et~al\mbox{.}}{2003}]%
        {mccormick2003easy}
\bibfield{author}{\bibinfo{person}{S~Thomas McCormick}, \bibinfo{person}{M~Rammohan Rao}, {and} \bibinfo{person}{Giovanni Rinaldi}.} \bibinfo{year}{2003}\natexlab{}.
\newblock \showarticletitle{Easy and difficult objective functions for max cut}.
\newblock \bibinfo{journal}{\emph{Mathematical programming}}  \bibinfo{volume}{94} (\bibinfo{year}{2003}), \bibinfo{pages}{459--466}.
\newblock


\bibitem[\protect\citeauthoryear{Mirka and Williamson}{Mirka and Williamson}{2023}]%
        {mirka2023experimental}
\bibfield{author}{\bibinfo{person}{Renee Mirka} {and} \bibinfo{person}{David~P Williamson}.} \bibinfo{year}{2023}\natexlab{}.
\newblock \showarticletitle{An Experimental Evaluation of Semidefinite Programming and Spectral Algorithms for Max Cut}.
\newblock \bibinfo{journal}{\emph{ACM Journal of Experimental Algorithmics}}  \bibinfo{volume}{28} (\bibinfo{year}{2023}), \bibinfo{pages}{1--18}.
\newblock


\bibitem[\protect\citeauthoryear{Nemhauser and Trotter~Jr}{Nemhauser and Trotter~Jr}{1975}]%
        {nemhauser1975vertex}
\bibfield{author}{\bibinfo{person}{George~L Nemhauser} {and} \bibinfo{person}{Leslie~E Trotter~Jr}.} \bibinfo{year}{1975}\natexlab{}.
\newblock \showarticletitle{Vertex packings: structural properties and algorithms}.
\newblock \bibinfo{journal}{\emph{Mathematical Programming}} \bibinfo{volume}{8}, \bibinfo{number}{1} (\bibinfo{year}{1975}), \bibinfo{pages}{232--248}.
\newblock


\bibitem[\protect\citeauthoryear{Oku, Tawada, Tanaka, and Togawa}{Oku et~al\mbox{.}}{2020}]%
        {oku2020reduce}
\bibfield{author}{\bibinfo{person}{Daisuke Oku}, \bibinfo{person}{Masashi Tawada}, \bibinfo{person}{Shu Tanaka}, {and} \bibinfo{person}{Nozomu Togawa}.} \bibinfo{year}{2020}\natexlab{}.
\newblock \showarticletitle{How to reduce the bit-width of an Ising model by adding auxiliary spins}.
\newblock \bibinfo{journal}{\emph{IEEE Trans. Comput.}} \bibinfo{volume}{71}, \bibinfo{number}{1} (\bibinfo{year}{2020}), \bibinfo{pages}{223--234}.
\newblock


\bibitem[\protect\citeauthoryear{Oshiyama and Ohzeki}{Oshiyama and Ohzeki}{2022}]%
        {oshiyama2022benchmark}
\bibfield{author}{\bibinfo{person}{Hiroki Oshiyama} {and} \bibinfo{person}{Masayuki Ohzeki}.} \bibinfo{year}{2022}\natexlab{}.
\newblock \showarticletitle{Benchmark of quantum-inspired heuristic solvers for quadratic unconstrained binary optimization}.
\newblock \bibinfo{journal}{\emph{Scientific reports}} \bibinfo{volume}{12}, \bibinfo{number}{1} (\bibinfo{year}{2022}), \bibinfo{pages}{2146}.
\newblock


\bibitem[\protect\citeauthoryear{Pardalos and Rodgers}{Pardalos and Rodgers}{1990}]%
        {pardalos1990computational}
\bibfield{author}{\bibinfo{person}{Panos~M Pardalos} {and} \bibinfo{person}{Gregory~P Rodgers}.} \bibinfo{year}{1990}\natexlab{}.
\newblock \showarticletitle{Computational aspects of a branch and bound algorithm for quadratic zero-one programming}.
\newblock \bibinfo{journal}{\emph{Computing}} \bibinfo{volume}{45}, \bibinfo{number}{2} (\bibinfo{year}{1990}), \bibinfo{pages}{131--144}.
\newblock


\bibitem[\protect\citeauthoryear{Powell}{Powell}{1994}]%
        {powell1994direct}
\bibfield{author}{\bibinfo{person}{Michael~JD Powell}.} \bibinfo{year}{1994}\natexlab{}.
\newblock \bibinfo{booktitle}{\emph{A direct search optimization method that models the objective and constraint functions by linear interpolation}}.
\newblock \bibinfo{publisher}{Springer}.
\newblock


\bibitem[\protect\citeauthoryear{Punnen}{Punnen}{2022}]%
        {punnen2022quadratic}
\bibfield{author}{\bibinfo{person}{Abraham~P Punnen}.} \bibinfo{year}{2022}\natexlab{}.
\newblock \bibinfo{booktitle}{\emph{The Quadratic Unconstrained Binary Optimization Problem: Theory, Algorithms, and Applications}}.
\newblock \bibinfo{publisher}{Springer Nature}.
\newblock


\bibitem[\protect\citeauthoryear{Rehfeldt, Koch, and Shinano}{Rehfeldt et~al\mbox{.}}{2023}]%
        {rehfeldt2023faster}
\bibfield{author}{\bibinfo{person}{Daniel Rehfeldt}, \bibinfo{person}{Thorsten Koch}, {and} \bibinfo{person}{Yuji Shinano}.} \bibinfo{year}{2023}\natexlab{}.
\newblock \showarticletitle{Faster exact solution of sparse MaxCut and QUBO problems}.
\newblock \bibinfo{journal}{\emph{Mathematical Programming Computation}} (\bibinfo{year}{2023}), \bibinfo{pages}{1--26}.
\newblock


\bibitem[\protect\citeauthoryear{Shih, Wu, and Kuo}{Shih et~al\mbox{.}}{1990}]%
        {shih1990unifying}
\bibfield{author}{\bibinfo{person}{Wei-Kuan Shih}, \bibinfo{person}{Sun Wu}, {and} \bibinfo{person}{Yue-Sun Kuo}.} \bibinfo{year}{1990}\natexlab{}.
\newblock \showarticletitle{Unifying maximum cut and minimum cut of a planar graph}.
\newblock \bibinfo{journal}{\emph{IEEE Trans. Comput.}} \bibinfo{volume}{39}, \bibinfo{number}{5} (\bibinfo{year}{1990}), \bibinfo{pages}{694--697}.
\newblock


\bibitem[\protect\citeauthoryear{Shor}{Shor}{2002}]%
        {shor2002introduction}
\bibfield{author}{\bibinfo{person}{Peter~W Shor}.} \bibinfo{year}{2002}\natexlab{}.
\newblock \showarticletitle{Introduction to quantum algorithms}. In \bibinfo{booktitle}{\emph{Proceedings of Symposia in Applied Mathematics}}, Vol.~\bibinfo{volume}{58}. \bibinfo{pages}{143--160}.
\newblock


\bibitem[\protect\citeauthoryear{Spielman and Srivastava}{Spielman and Srivastava}{2008}]%
        {spielman2008graph}
\bibfield{author}{\bibinfo{person}{Daniel~A Spielman} {and} \bibinfo{person}{Nikhil Srivastava}.} \bibinfo{year}{2008}\natexlab{}.
\newblock \showarticletitle{Graph sparsification by effective resistances}. In \bibinfo{booktitle}{\emph{Proceedings of the fortieth annual ACM symposium on Theory of computing}}. \bibinfo{pages}{563--568}.
\newblock


\end{thebibliography}
